\documentclass{elsart}
\usepackage{amsfonts}

\newcommand{\eq}[1]{\textrm{(#1)}}
\newcommand{\ii}{\mathrm{i}}
\newcommand{\longto}{\mathop{\longrightarrow}\limits}
\newcommand{\opname}[1]{\mathop{\mathrm{#1}}\nolimits}
\newcommand{\pd}[2]{\frac{\partial#1}{\partial#2}}
\newcommand{\sepword}[1]{\qquad\mbox{#1}\quad}
\newcommand{\set}[1]{\{\,#1\,\}}
\DeclareMathSymbol{\toto}{\mathrel}{AMSa}{"13}

\begin{document}

\begin{frontmatter}


\title{Connes' tangent groupoid and strict quantization}

\author[FTUZ]{Jos\'e F. Cari\~nena},
\author[FTUZ]{Jes\'us Clemente-Gallardo},
\author[FTUZ]{Eduardo Follana},
\author[FTUZ,FUCR]{Jos\'e M. Gracia-Bond\'{\i}a},
\author[FTUZ]{Alejandro Rivero}, and
\author[MUCR]{Joseph C. V\'arilly}

\address[FTUZ]{Departamento de F\'{\i}sica Te\'orica,
               Universidad de Zaragoza, 50009 Zaragoza, Spain}
\address[FUCR]{Departamento de F\'{\i}sica,
               Universidad de Costa Rica, 2060 San Pedro, Costa Rica}
\address[MUCR]{Departamento de Matem\'aticas,
               Universidad de Costa Rica, 2060 San Pedro, Costa Rica}

\begin{abstract}
We address one of the open problems in quantization theory
recently listed by Rieffel. By developing in detail Connes' tangent
groupoid principle and using previous work by Landsman, we show how to
construct a strict flabby quantization, which is moreover an
asymptotic morphism and satisfies the reality and traciality
constraints, on any oriented Riemannian manifold. That construction
generalizes the standard Moyal rule. The paper can be considered as an
introduction to quantization theory from Connes' point of view.
\end{abstract}

\end{frontmatter}

\section{Introduction}
\label{sc:one}

\subsection{Motivation}

Recently, Rieffel has published a list of open
problems~\cite{RieffelProgram} in quantization. The main aim of this
paper is to address Question~20 of Rieffel's list: ``$\dots\,$in what
ways can a suitable Riemannian metric on a manifold $M$ be used to
obtain a strict deformation quantization of $T^*M$?'' We do that by
giving proofs for (a slightly improved variant of) a construction
sketched by Connes in Section~II.5 of his book~\cite{Book} on
noncommutative geometry, and elsewhere~\cite{ConnesLH}. The paper can
be considered as an introduction to the subject of quantization from
Connes' point of view and thus serves too a pedagogical purpose.

Rieffel's requirements are stronger that those of formal deformation
theory, extant for any Poisson manifold~\cite{Kontsevich}. To our
mind, however, the fact that $T^*M$, for $M$ Riemannian, possesses a
strict quantization, was indeed proved by Landsman in the
path-breaking paper~\cite{Landsman}. Nevertheless, the noncommutative
geometry approach presents several advantages, not least that the
$C^*$-theoretical aspects come in naturally. The procedure was
suggested by Landsman himself, even before~\cite{Book} was in print,
at the end of his paper. For the sake of simplicity, we deal here with
the nonequivariant case only.

The plan of the article is as follows. In this first Section, after
introducing groupoids and the tangent groupoid construction, we give
an elementary discussion of Connes' recipe for quantization. We find
the intersection, for the case $M = \Rset^n$, of the family of what
could be termed ``Connes' quantization rules'' with the ordinary
quantum-mechanical ordering prescriptions. Groupoids are here regarded
set-theoretically, questions of smoothness being deferred to
Section~\ref{sc:two}. We show that Moyal's quantization rule belongs
to the collection of Connes' quantization rules and in fact is singled
out by natural conditions in strict deformation theory.

In Section~\ref{sc:two} we spell out our variant of Connes' tangent
groupoid construction in full detail. The heart of the matter is the
continuity of the groupoid product, for which we give two different
proofs.

In Section~\ref{sc:three}, by means of the mathematical apparatus of
Section~\ref{sc:two}, we restate Landsman's partial answer to
Rieffel's question. In particular, we rework the existence proof for a
strict quantization \textit{of the Moyal type} (in the sense of being
both real and tracial) on Riemannian manifolds. Moreover, using a
strong form of the tubular neighbourhood theorem, we show that there
exists what Rieffel calls a \textit{flabby}
quantization~\cite{RieffelProgram}. The paper concludes with a
discussion on the $C^*$-algebraic aspect of the tangent groupoid
construction and its relation to the index theorem.

\subsection{Basic facts on groupoids}

The most economical way to think of a groupoid is as a pair of sets
$G^0 \subset G$, and to regard elements of $G$ as arrows and elements
of $G^0$ as nodes.

\begin{defn}
A \textit{groupoid} $G \toto G^0$ is a small category in which every
morphism has an inverse. Its set of objects is~$G^0$, its set of
morphisms is~$G$.
\end{defn}

A group is of course a groupoid with a single object. The gist of the
definition is conveyed by the following example.

\paragraph{Partial isometries}
Consider a complex Hilbert space $\mathcal{H}$. The collection of
\textit{unitary} arrows between closed subspaces of~$\mathcal{H}$
obviously defines a groupoid, for which $G$ is the set of partial
isometries $\set{w \in \mathcal{L}(\mathcal{H}) : ww^*w = w}$ and
$G^0$ is the set of orthogonal projectors in
$\mathcal{L}(\mathcal{H})$. Recall that, given $w$, we can write
\[
\mathcal{H} = \ker w \oplus \opname{im} w^*
 = \ker w^* \oplus \opname{im} w.
\]
Hence $w^*w$ is the orthogonal projector with range $\opname{im}w^*$,
while $ww^*$ is the orthogonal projector with range $\opname{im}w$;
of course, $w^*$ is the inverse of~$w$. We naturally identify two maps
$r$, $s$ (respectively ``range'' and ``source'') from $G$ to $G^0$:
$r(w) := ww^*$ and $s(w) := w^*w$. Also, it is natural to consider
\[
G^{(2)} := \set{(u,v) \in G \times G : u^*u = vv^*}.
\]
This defining equation is a sufficient condition for the operator
product $uv$ to be a partial isometry, as follows from the simple
calculation
\[
uv(uv)^*uv = uvv^*u^*uv = uvv^*v = uv.
\]

The example motivates a more cumbersome restatement, that includes all
the practical elements of the definition.

\begin{defn}
\label{df:gpoid}
A groupoid $G \toto G^0$ consists of: a set $G$, a set $G^0$ of
``units'' with an inclusion $G^0 \hookrightarrow G$, two maps
$r,s : G \to G^0$, and a composition law $G^{(2)} \to G$ with domain
\[
G^{(2)} := \set{(g,h) : s(g) = r(h)} \subseteq G \times G,
\]
subject to the following rules:
\begin{enumerate}
\item if $g \in G^0$ then $r(g) = s(g) = g$;
\item $r(g)g = g = gs(g)$;
\item each $g \in G$ has an ``inverse'' $g^{-1}$, satisfying
      $gg^{-1} = r(g)$ and $g^{-1}g = s(g)$;
\item $r(gh) = r(g)$ and $s(gh) = s(h)$ if $(g,h) \in G^{(2)}$;
\item $(gh)k = g(hk)$ if $(g,h) \in G^{(2)}$ and $(gh,k) \in G^{(2)}$.
\end{enumerate}
\end{defn}

The examples of groupoids we shall use have a differential geometric
flavour instead.

\paragraph{A vector bundle $E \longto^\pi M$}
Here $G = E$ is the total space, $G^0 = M$ is the base space,
$r = s = \pi$ so that $G^{(2)} = \biguplus_{x\in M} E_x \times E_x$
(the total space of the Whitney sum $E \oplus E$), and the composition
law is \textit{fibrewise addition}.

\paragraph{The double groupoid of a set}
Given a set $M$, take $G = M \times M$ and $G^0 = M$, included in
$M \times M$ as the diagonal subset
$\Delta(M) := \set{(x,x) : x \in M}$. Define $r(x,y) := x$,
$s(x,y) := y$. Then $(x,y)^{-1} = (y,x)$ and the composition law is
\[
(x,y) \cdot (y,z) = (x,z).
\]
We shall generally call $G^0$ the \textit{diagonal} of $G$.

\subsection{Connes' tangent groupoids}

We recall some basic concepts of differential geometry, particularly
sprays and normal bundles, that we shall later use.

Given a symmetric linear connection on a differentiable manifold $M$,
one can define a vector field $\Gamma$ on~$TM$ whose value at
$v \in TM$ is its horizontal lift to~$T_v(TM)$. This vector field is
called the geodesic spray of the connection and its integral curves
are just the natural lift of geodesics in~$M$. The geodesic spray is a
second-order differential equation vector field satisfying an
additional condition of degree-one homogeneity which corresponds to
the affine reparametrization property of geodesics. More generally, a
second-order differential equation vector field is said to be a spray
if the set of its integral curves is invariant under any affine
reparametrization: these curves are called geodesics of the spray.
Given a spray $\Gamma$, there is a symmetric connection whose geodesic
spray is~$\Gamma$; and conversely, the connection is fully determined
by its geodesic spray, that can also be used to construct the
exponential map $\exp\colon T_xM \to M$. A Riemannian structure on~$M$
determines one symmetric linear connection, the Levi-Civita
connection, and consequently a Riemannian spray.

Now let $Y^0$ be a submanifold of a manifold $Y$. The normal bundle
$\mathcal{N}_{Y^0}^Y$ to $Y^0$ in $Y$ is defined as the vector bundle
$\mathcal{N}_{Y^0}^Y := TY/TY^0$, where the notation means that its
base is $Y^0$ and its fibre is given by the equivalence classes of the
elements of the tangent bundle $TY$ under the relation:
$X_1 \sim X_2$, for $X_1,X_2 \in T_q Y$ with $q \in Y^0$, if and only
if $X_1 = X_2 + V$ for some $V \in T_q Y^0$. The usual way to work
with such a structure is to choose a representative in each class,
thereby forming a complementary bundle to $TY^0$ in $TY$ restricted
to~$Y^0$. There is no canonical choice for the latter, in general.
When a Riemannian metric is provided on~$Y$, there is a natural
definition of the complementary bundle as the orthogonal complement of
the tangent space $T_qY^0$ in $T_qY$ (although other choices may be
convenient, even in the Riemannian case). Once we have chosen a
suitable representative of each class, the bundle
$\mathcal{N}_{Y^0}^Y$ becomes a subbundle of $TY$ and we can consider
the exponential map $\exp$ of $TY$ restricted
to~$\mathcal{N}_{Y^0}^Y$.

We recall also that a \textit{tubular neighbourhood} of $Y^0$ in $Y$
is a vector bundle $E \to Y^0$, an open neighbourhood $Z$ of its zero
section and a diffeomorphism of $Z$ onto an open set (the tube)
$U \subset Y$ containing $Y^0$, which restricts over the zero section
to the inclusion of $Y^0$ in $Y$~\cite{Lang}. We say that the tubular
neighbourhood is \textit{total} when $Z = E$. The main theorem in this
context establishes that, given a spray on~$Y$, one can always
construct a tubular neighbourhood by making use of the corresponding
exponential map. When the spray is associated to a Riemannian metric,
one can always have a total tubular neighbourhood, because a Euclidean
bundle is compressible, i.e., isomorphic as a fibre bundle to an open
neighbourhood of the zero section: see~\cite{Lang}.

For suitable $\hbar_0 > 0$, one can then define the normal cone
deformation (a sort of blowup in the differentiable category) of the
pair $(Y,Y^0)$, denoted $\mathcal{M}_{Y^0}^Y$, by gluing together
$Y \times (0,\hbar_0]$ with $\mathcal{N}_{Y^0}^Y$ as a boundary with
the help of the tubular neighbourhood construction~\cite{HilsumS}. The
construction is particularly interesting when $(Y,Y^0)$ is a groupoid,
in that it gives a ``normal groupoid'' with diagonal
$Y^0 \times [0,\hbar_0]$.

We consider a particular case of this construction, the tangent
groupoid. Let $M$ now be an orientable Riemannian manifold. Denote by
$\mathcal{N}^\Delta$, rather than $\mathcal{N}_M^{M\times M}$, the
normal bundle associated to the diagonal embedding
$\Delta\colon M \to M \times M$. We can identify $TM \oplus TM$ with
the restriction of $T(M \times M)$ to $\Delta(M)$, and the tangent
bundle over $M$ is identified to
\[
\set{(\Delta(q), X_q, X_q) : (q,X_q) \in TM};
\]
thereby the normal bundle $\mathcal{N}^\Delta$ to~$M$ in $M \times M$
can \textit{a priori} be identified with
\begin{equation}
\set{(\Delta(q), \varphi_{1q}X_q, \varphi_{2q}X_q) : (q,X_q) \in TM},
\label{eq:normal-bundle}
\end{equation}
where $\varphi_1,\varphi_2 \in \opname{End}(TM)$ are any two bundle
endomorphisms (i.e., continuous vector bundle maps from $TM$ into
itself) such that the linear map $\varphi_{1q} - \varphi_{2q}$
on~$T_qM$ is invertible for all $q \in M$; this we write as
$\varphi_1 - \varphi_2 \in GL(TM)$. We shall assume, for definiteness,
that each $\varphi_{1q} - \varphi_{2q}$ is homotopic to the identity:
in other words, that there is an isomorphism of oriented vector
bundles between $TM$ and~$\mathcal{N}^\Delta$.

The tangent groupoid $G_M \toto G^0_M$, according to Connes, is
essentially the normal groupoid $\mathcal{M}_M^{M\times M}$ modulo
that isomorphism~\cite{Book}. That is to say, we think of the disjoint
union $G_M = G_1 \uplus G_2$ of two groupoids
\[
G_1 := M \times M \times (0,\hbar_0],  \qquad  G_2 := TM,
\]
where $G_1$ is the disjoint union of copies
$M \times M \times \{\hbar\}$ of the double groupoid of~$M$,
parametrized by $0 < \hbar \leq \hbar_0$. The compositions are
\[
\begin{array}{rcl}
(x,y,\hbar) \cdot (y,z,\hbar) & = & (x,z,\hbar)     \\
   (q, X_1) \cdot (q, X_2)    & = & (q, X_1 + X_2)
\end{array}
\qquad
\begin{array}{l}
\mbox{with \ } 0 < \hbar \leq \hbar_0, \\
\mbox{if \ }   X_1, X_2 \in T_qM.
\end{array}
\]
The topology on $G_M$ is such that, if $(x_n, y_n, \hbar_n)$ is a
sequence of elements of $G_1$ with $\hbar_n \downarrow 0$, then it
converges to a tangent vector $(q,X)$ iff
\[
x_n \to q, \quad  y_n \to q
\]
and
\begin{equation}
\frac{x_n - y_n}{\hbar_n} \to X.
\label{eq:bdry-tgt}
\end{equation}
The last condition, although formulated in a local chart, clearly
makes intrinsic sense, as $x_n - y_n$ and $X$ get multiplied by the
same factor, to first order, under a change of charts. It forces us to
demand that $\varphi_{1q} - \varphi_{2q}$ be equal to the identity in
\eq{\ref{eq:normal-bundle}}; this will be seen later to be necessary
also for linear symplectic invariance of the associated quantization
recipe in the flat case. Actually Connes picks, to coordinatize the
tangent groupoid, the particular identification given by
$\varphi_{1q}X_q = 0$, $\varphi_{2q}X_q = - X_q$.

\subsection{Connes' quantization rule on $\Rset^n$}

Consider then the case $M = \Rset^n$. The groupoid $G_{\Rset^n}$ is
certainly a manifold with boundary $T\Rset^n$. The neighbourhood of
points in the boundary piece is described as follows: for
$(q,X) \in T\Rset^n$, define
\[
(x, y, \hbar) := \Phi_\hbar(q, X, \hbar)
 = (q + \hbar\varphi_{1}X, q + \hbar\varphi_{2}X, \hbar),
\]
with $\varphi_{1}X - \varphi_{2}X = X$. (It is natural in this context
to take $\varphi_{1q}$, $\varphi_{2q}$ to be constant maps.) In fact,
$\Phi$ is a linear isomorphism of $T\Rset^n \times [0,1]$ onto the
groupoid.

Now comes Connes' quantization recipe: a function on $G_M$ is first of
all a pair of functions on $G_1$ and $G_2$ respectively. The first one
is essentially the kernel for an operator on $L^2(M)$, the second the
Fourier transform of a function on the cotangent bundle $T^*M$, i.e.,
the classical phase space. The condition that both match seamlessly to
give a regular function on $G_M$ is precisely construed as the
quantization rule!

Let $a(q,p)$ be a function on $T^*\Rset^n$. Its inverse Fourier
transform in the second variable gives us a function on~$T\Rset^n$:
\[
\mathcal{F}^{-1} a(q,X) = \int_{\Rset^n} e^{\ii Xp} a(q,p) \,\d p.
\]
To the function~$a$, Connes' prescription associates then the
following family of kernels:
\[
k_a(x,y;\hbar) := \mathcal{F}^{-1} a(\Phi_\hbar^{-1}(x,y,\hbar)).
\]
The map $a \mapsto k_a$ is thus linear. We get the dequantization rule
by Fourier inversion:
\begin{equation}
a(q,p) = \frac{1}{(2\pi\hbar)^n} \int_{\Rset^n}
           k_a(q + \hbar\varphi_{1}X, q + \hbar\varphi_{2}X; \hbar)
            \, e^{-\ii pX} \,\d X.
\label{eq:Fourier-inv}
\end{equation}

We may well ask where Connes' rule stands with respect to the usual
quantization rules. We find it more convenient to argue from the
dequantization formula. Any linear dequantization rule can be
expressed in the form
\[
a(q,p)
 = \int_{\Rset^{2n}} K(q,p,x,y;\hbar)\, k_a(x,y;\hbar) \,\d x\,\d y,
\]
for a suitable distributional kernel $K$. To select useful rules, one
seeks to impose reasonable conditions. The simplest is obviously
equivariance under translations in $T^*\Rset^n$, corresponding
physically to Galilei invariance. Equivariance under translation of
the spatial coordinates amounts immediately to the condition
\begin{equation}
\biggl( \pd{}{q} + \pd{}{x} + \pd{}{y} \biggr) K = 0,
\label{eq:derK-zero}
\end{equation}
i.e., $K$ depends only on the combinations $x - y$ and
$q - \half(x + y)$. A similar argument in the dual space shows that
equivariance under translation of the momenta amounts to
\begin{equation}
K(q, p + p', x,y; \hbar) = e^{-\ii p'(x-y)/\hbar} K(q,p,x,y; \hbar).
\label{eq:phase-K}
\end{equation}
As beautifully discussed in~\cite{KrugerP}, both conditions together
call for the role of Weyl operators in quantization. In fact, it is
clear that the most general kernel satisfying \eq{\ref{eq:derK-zero}}
and \eq{\ref{eq:phase-K}} is of the form
\[
K(q, p, x, y; \hbar)
= \frac{1}{(2\pi\hbar)^{2n}} e^{-\ii p(x-y)/\hbar} \int_{\Rset^n}
  f(\theta, x-y) \exp(-\ii\theta(q - \half(x+y))/\hbar) \,\d\theta,
\]
where $f$ should have no zeros, on account of invertibility; that can
be rewritten as
\[
K(q, p, x, y; \hbar) = \frac{1}{(2\pi\hbar)^{2n}} \int_{\Rset^{2n}}
 f(\theta, \tau) e^{-\ii(\theta q + \tau p)/\hbar}
  \langle y \mathbin{\vert}
    \exp(-\ii(\theta Q + \tau P)/\hbar)x \rangle \,\d\theta \,\d\tau,
\]
where $Q$, $P$ denote the usual quantum observables for position and
momentum, respectively, and $\exp(-\ii(\theta Q + \tau P)/\hbar)$ are
the Weyl operators.

Further requirements on the quantization rule lead to a sharpening in
the determination of~$f$. For instance, in the terminology
of~\cite{Iapetus}, the requirement of \textit{semitraciality} (that
dequantization of the kernel corresponding to the identity operator be
the function~$1$) leads to $f(0,0) = 1$; of \textit{reality} (that the
kernels for real classical observables correspond to selfadjoint
operators (quantum observables according to Quantum Mechanics
postulates) forces $f(\theta,\tau) = f^*(-\theta, -\tau)$, and so on.
Moyal's rule~\cite{Moyal}, that is, the paradigmatic example of strict
deformation quantization, corresponds to taking $f = 1$.

A very important requirement, still in the terminology
of~\cite{Iapetus}, is \textit{traciality}; that is, the coincidence
between the classical and quantum averages for the product of
observables. Mathematically it is equivalent to unitarity of $K$,
which demands $|f|^2 = 1$. Taken together, traciality and reality
allow to rewrite quantum mechanics as a statistical theory in
classical phase space; and mathematically, they practically force the
Moyal rule. Traciality is further discussed in Section~\ref{sc:three},
in the general context of Riemannian manifolds.

On the other hand, Connes' dequantization rule
\eq{\ref{eq:Fourier-inv}} gives
$$
K(q, p, x, y; \hbar) = \frac{1}{(2\pi\hbar)^n} \int_{\Rset^n}
 \delta(x-q-\hbar\varphi_1 X) \delta(y-x+\hbar X) \exp(-\ii pX)\,\d X.
$$
This integral can be rewritten as
$$
(2\pi\hbar)^{-n} \,\delta(x-q-\varphi_1(x-y))\, e^{-\ii p(x-y)/\hbar},
$$
and further transformed into
$$
\frac{1}{(2\pi\hbar)^{2n}} e^{-\ii p(x-y)/\hbar} \int_{\Rset^n}
 \exp(\ii\theta(\half(x+y)-q+(\half-\varphi_1)(x-y))/\hbar)\,\d\theta,
$$
which is of the previously given form, with
$$
f(\theta,\tau) = \exp(\ii\theta(\half-\varphi_1)\tau/\hbar).
$$
We conclude that Connes' quantization rule, in the flat space case, is
determined up to an arbitrary linear transformation of $\Rset^n$.

Had we tried to keep the most general formula
$$
K(q, p, x, y; \hbar) = \frac{1}{(2\pi\hbar)^n} \int_{\Rset^n}
 \delta(x-q-\hbar\varphi_{1q}X)
 \delta(y-x+\hbar(\varphi_{1q}-\varphi_{2q})X) \exp(-\ii pX) \,\d X,
$$
we would have had two difficulties: equivariance under translation of
the position would be lost since $\varphi_{1q} \neq \varphi_{1q'}$ in
general, and equivariance under translation of the momentum would be
lost since $\varphi_{1q} - \varphi_{2q} \neq \opname{id}_{T_q\Rset^n}$
in general; this bears on the physical meaning of Connes' limit
condition~\eq{\ref{eq:bdry-tgt}}.

It would be tempting, also in view of \eq{\ref{eq:bdry-tgt}}, to add
the time variable to the mathematical apparatus of this paper, in the
spirit of Feynman's formalism, and to study the interchangeability of
the limits $t \downarrow 0$ and $\hbar \downarrow 0$.

Happily, Moyal's rule is included among Connes' rules: it follows from
the most natural choice $\varphi_1 = \half$, $\varphi_2 = - \half$
(see the discussion in~\cite{Portia}). We arrive at the following
conclusion.

\begin{prop}
Moyal quantization rule is the only \textbf{real} quantization of the
Connes type, in the case $M = \Rset^n$.   
\end{prop}

Connes' own choice in~\cite{Book} is $\varphi_1 = 0$,
$\varphi_2 = -1$; this corresponds to the ``standard'' ordering
prescription, in which the quantization of $q^np^m$ is $Q^nP^m$;
whereas the choice $\varphi_1 = 1$, $\varphi_2 = 0$ leads to the
``antistandard'' ordering, in which the quantization of $p^mq^n$ is
$P^mQ^n$. Note that all of Connes' prescriptions are tracial. In order
to obtain the ``normal'' and ``antinormal'' prescriptions of use in
field theory, which are real but not tracial, one would have to
complexify Connes' construction; we shall not go into that.

Following Landsman, we will strengthen the definition of strict
quantization in order to remain in the real context (see
Section~\ref{sc:three}). Our elementary discussion in this subsection
leads to presume that real strict quantization can be done in the
framework of tangent groupoids; this is tantamount to the
generalization of the Moyal rule to arbitrary manifolds endowed with
sprays. All the pertinent differential geometric constructions for
that purpose are taken up in the next section.

\section{The tangent groupoid construction: the Moyal version}
\label{sc:two}

\begin{defn}
A \textbf{smooth groupoid} is a groupoid $(G,G^0)$ together with
differentiable structures on $G$ and $G^0$ such that the maps $r$ and
$s$ are submersions, and the inclusion map $G^0 \hookrightarrow G$ is
smooth as well as the product $G^{(2)} \to G$.
\end{defn}

Note the dimension count: if $\dim G = n$, $\dim G^0 = m$, then
$\dim G^{(2)} = 2n - m$.

The definition of smooth groupoid is taken from~\cite{Book}. Now we
establish his Proposition II.5.4, there left unproven by Connes.

\begin{prop}
\label{pr:smooth-gpoid}
The tangent groupoid $G_M$ to a smooth manifold $M$ is a smooth
groupoid.
\end{prop}

\subsection{$G_M$ as a manifold with boundary}

The groupoid $G_1 = M \times M \times (0,\hbar_0]$, that will be the
interior of~$G_M$ (plus a trivial ``outer'' boundary, which we neglect
to mention in our following arguments) is given the usual product
manifold structure. To complete the definition of the manifold
structure of $G$, consider the isomorphism from $TM$ to $N^\Delta$
given by $(q,X_q) \mapsto (q, \half X_q, -\half X_q)$. Consider also
the product manifold $TM \times [0,\hbar_0]$. We choose a spray on $M$
---this is provided, for instance, by a choice of Riemannian metric
on~$M$--- and define a map
\[
TM \times [0,\hbar_0] \supset U \longto^{\Phi} G_M,
\]
where $U$ is open in $TM \times [0,\hbar_0]$ and includes
$TM \times \{0\}$, as follows:
\[
\Phi(q, X_q, \hbar)
 := (\exp_q(\hbar X_q/2), \exp_q(-\hbar X_q/2), \hbar)
 \sepword{for} \hbar > 0,
\]
and
\[
\Phi(q, X_q, 0) := (q, X_q)  \sepword{for} \hbar = 0.
\]
Here $\exp$ denotes the exponential map associated to the spray: we
know that, for a fixed $\hbar$, the exponential map defines a
diffeomorphism of an open neighbourhood $V$ of the zero section of
$TM$ onto an open neighbourhood of the diagonal in $M \times M$; and
we decide that a point $(q,X_q,\hbar)$ is in $U$ if
$\Phi(q,X_q,\hbar)$ is contained in $V$.

Therefore, both the existence and, for a suitable choice of $U$, the
bijectivity of the map $\Phi$ follow from the tubular neighbourhood
theorem. As $U$ is an open (sub)manifold with boundary, we can carry
the structure of manifold with boundary to $G_M$, obtaining that $TM$
is the boundary of the groupoid $G_M$ in the topology associated to
that structure. The diagonal is obviously $M \times [0,\hbar_0]$.

We remind the reader that Connes uses instead the chart given by
\[
(q, X_q, \hbar) \mapsto (q, \exp_q(-\hbar X_q), \hbar)
 \sepword{for} \hbar > 0,
\]
and
\[
(q, X_q, 0) \mapsto (q, X_q)  \sepword{for} \hbar = 0.
\]

{}From now on, a particular Riemannian structure is assumed chosen,
and when convenient we shall also assume, as we then may, that $U$ is
the whole of $TM \times [0,\hbar_0]$. To continue the proof of
smoothness of $G_M$, we need to check that the various mappings have
the required smoothness properties. The basic idea is to pull all maps
on $G_M$ back to $TM \times [0,\hbar_0]$ and prove smoothness there.

Consider first the inclusion map $i\colon G^0_M \to G_M$. It is
obvious that the restriction of $i$ to $M \times (0,\hbar_0]$ is
smooth in its domain. If we now consider its restriction to
$i^{-1}(\Phi(U))$ and compose it with $\Phi^{-1}$, we obtain a map
which can be written as:
\[
\Phi^{-1} \circ i(x,\hbar) = \Phi^{-1}(x, x, \hbar) = (x, 0_x, \hbar)
 \sepword{for} \hbar > 0,
\]
and
\[
\Phi^{-1} \circ i(x,0) = \Phi^{-1}(x, 0_x, 0) = (x, 0_x, 0)
  \sepword{for} \hbar = 0.
\]
This map is obviously smooth in its domain and, as $\Phi$ is a
diffeomorphism, $i$~is smooth in its domain.

We shall consider now the range and source maps. The smoothness of
both maps when restricted to $M \times M \times (0,\hbar]$ is again
obvious, as is the fact that they are of maximum rank (hence
submersions when restricted to this domain). The composition of $\Phi$
with the restriction of $r$ to $\Phi(U)$ is expressed as:
\[
r \circ \Phi (q, X_q, \hbar)
 = r(\exp_q(\hbar X_q/2), \exp_p(-\hbar X_p/2), \hbar)
 = (\exp_q(\hbar X_q/2),\hbar)
\]
for $\hbar > 0$ and
\[
r \circ \Phi (q, X_q, 0) = (q, 0)
\]
for $\hbar = 0$. Again this map is smooth and of maximum rank in its
domain, so that $r$ must also be a submersion. The corresponding proof
for the source map is analogous.

\subsection{The geometrical structure of
            $G^{(2)}_M \subset G_M \times G_M$}

To define a differentiable structure for $G^{(2)}_M$, we proceed as in
the previous case, by defining a bijection between an open set in
$G^{(2)}_M$ and an open set in a manifold with boundary, and
transporting the differential structure via the bijection.

What kind of manifold is $G^{(2)}_M$? The product has to be a smooth
mapping between two manifolds with boundary, mapping the boundary on
the boundary and the interior on the interior. Thus, from the
definition of the product it is clear that the boundary of the
manifold should be the Whitney sum $TM \oplus TM$, arising as the
pullback with respect to the diagonal injection of $M$ in $M \times M$
of the product bundle $TM \times TM$.

Points in the interior of the manifold are pairs of the form:
$((x,y,\hbar),(y,z,\hbar))$. Therefore, we need to define a
differential structure in $G^{(2)}_M$ in such a way that it becomes a
manifold with boundary $TM \oplus TM$, its interior being
diffeomorphic to $M \times M \times M \times (0,\hbar_0]$.

Consider now $TM \times TM \times [0,\hbar_0]$. This is a manifold
with boundary, with a natural differential structure. Let
$(q',q,X_{q'},Y_q,\hbar)$ be a point in this manifold; we use our
construction for $G_M$ on each $TM$ separately, i.e., we set the
bijection:
\[
(q', q, X_{q'}, Y_q, \hbar) \leftrightarrow
 \bigl\{ (e_{q'}^{\hbar X_{q'}/2}, e_{q'}^{-\hbar X_{q'}/2}, \hbar),
   (e_q^{\hbar Y_q/2}, e_q^{-\hbar Y_q/2}, \hbar) \bigr\}
\]
for $\hbar > 0$ ---with an obvious notation for $\exp$. Those points
that interest us result just from imposing that
$e_{q'}^{-\hbar X_{q'}} = e_q^{\hbar Y_q}$. We shall see that this
constraint defines a regular submanifold of
$TM \times TM \times [0,\hbar_0]$. For that, for a fixed value of
$\hbar$, consider the sequence of maps
\[
\Psi : TM \times TM \longto^\Phi M \times M \times M \times M
 \longto^{\pi_{23}} M \times M \longto^{\Phi^{-1}} TM \to \Rset^n
\]
given by
\begin{eqnarray*}
(q', q, X_{q'}, Y_q)
& \mapsto & (e_{q'}^{\hbar X_{q'}/2}, e_{q'}^{-\hbar X_{q'}/2},
             e_q^{\hbar Y_q/2}, e_q^{-\hbar Y_q/2})
\\
& \mapsto & (e_{q'}^{-\hbar X_{q'}/2}, e_q^{\hbar Y_q/2})
             \mapsto (r, X_r) \mapsto X_r,
\end{eqnarray*}
where the $(r,X_r)$, that depend on $\hbar$, are found so $\Phi(r,X_r)
 = \bigl( e_{q'}^{-\hbar X_{q'}/2}, e_q^{\hbar Y_q/2} \bigr)$. This
composition defines a differentiable mapping of constant rank (it is
composed of two bijections and two projections onto factors of a
product). We then extend $\Psi$ to a map from
$TM \times TM \times [0,\hbar_0]$ to $\Rset^n \times [0,\hbar_0]$ and
have then that $\Psi^{-1}(\{0\} \times [0,\hbar_0])$ defines a regular
submanifold $S$ of $TM \times TM \times [0,\hbar_0]$, whose
differentiable structure we use to define the structure of manifold
with boundary on $G^{(2)}$. [It is perhaps not entirely clear that the
boundary of $S$ is $TM \oplus TM$, as we wish. But remember that $r$
depends on~$\hbar$ through the above manipulations. If we have a
sequence
\[
(q'_n, q_n, X_{q'_n}, Y_{q_n}, \hbar_n)
\]
in $S$, with $\hbar_n \downarrow 0$, its limit is a point $(s,s,X,Y)$
in $TM \times TM$ where $s = \lim q'_n = \lim q_n = \lim r(\hbar_n)$.]

\subsection{Continuity of the product}

A very simple argument with Riemannian flavour allows one to prove at
least continuity of the product operation in the tangent groupoid.
Consider a sequence in $G^{(2)}$:
\[
\{(x_n,y_n,\hbar_n), (y_n,z_n,\hbar_n)\}
 = \Phi(q'_n,q_n,X_n,Y_n,\hbar_n)
\]
with limit on the boundary. We need to check that the limit of the
products coincides with the product of the limits in $G_M$. For
elements close enough to the boundary, we know that $(s(n), Z_s(n))$
exists, so
\[
x_n = e_{s(n)}^{\hbar Z_s(n)/2},  \quad
z_n = e_{s(n)}^{-\hbar Z_s(n)/2}.
\]
Assume that
\[
q'_n \to s, \quad  q_n \to s, \quad  X_n \to A, \quad  Y_n \to B
 \qquad\mbox{as}\enspace  \hbar \downarrow 0.
\]
Now, we have
\[
\bigl( e_{s(n)}^{\hbar Z_s(n)/2}, e_{s(n)}^{-\hbar Z_s(n)/2} \bigr)
 = \bigl( e_{q'_n}^{\hbar X_n/2}, e_{q_n}^{-\hbar Y_n/2} \bigr)
 \to (s,s)  \qquad\mbox{as}\enspace  \hbar \downarrow 0,
\]
hence
\[
\lim s(n) = s.
\]
It remains to show that $\lim Z_s(n) = A + B$. This will follow if we
prove $Z_n = X_n + Y_n + o(\hbar_n)$. But that follows from
consideration of the small triangle with vertices $x_n$, $y_n$, $z_n$,
formed by geodesics through $q'_n$, $q_n$, $s(n)$ with directions
$\pm X_n$, $\pm Y_n$, $\pm Z_s(n)$.

One sees that $2(X_n + Y_n - Z_s(n))$ is approximately a circuit
around this triangle; by the Gauss--Bonnet theorem, we conclude that
\[
X_n + Y_n - Z_s(n) \approx O(\hbar^2).
\]

\subsection{A functorial proof of smoothness}

\begin{lem}
Let there be given two closed submanifolds $X_0 \hookrightarrow X$ and
$Y_0 \hookrightarrow Y$ and a smooth mapping $f\colon X \to Y$ which
satisfies $f(X_0) \subset Y_0$. Then the induced mapping
$\tilde f \colon \mathcal{M}_{X^0}^X \to \mathcal{M}_{Y^0}^Y$ between
the corresponding normal cone deformations is also smooth.
\end{lem}

This is Lemma~2.1 in~\cite{MonthubertAmiPierrot}, where no proof is
offered. We give some details of the lemma and then of its
application. First of all, if $f$ is a smooth mapping from $X$ to~$Y$
such that $f(X_0) \subset Y_0$, the image under $f_*$ of the tangent
bundle to $X_0$ is a subbundle of $TY_0$, implying the existence of an
induced mapping between the respective normal bundles, that we shall
continue to call~$f_*$. Now we define $\tilde f$ by
\[
\tilde f(x, \hbar) := (f(x), \hbar)  \sepword{for} \hbar > 0,
\]
and
\[
\tilde f(a, X_a) := (f(a), f_*(X_a))
\]
for $(a, X_a)$ an element of the normal bundle and $\hbar = 0$. Now,
differentiability of $\tilde f$ follows from the limit
\[
\lim_{\hbar\downarrow 0} \hbar^{-1}
 \exp_{f(a)}^{-1} \bigl[ f(\exp_a(\hbar X_a)) \bigr] = f_*(X_a).
\]

The application to the groupoid operations in our context is plain. We
discuss the product operation, the only relatively tricky one. Let
$m \colon G^{(2)} \to G$ be the groupoid multiplication. Now
$G^{(2)} \cap (G^0 \times G^0) = \Delta(G^0)$ and
$m(G^{(2)} \cap (G^0 \times G^0)) = G^0$; indeed, if
$(u,v) \in G^{(2)} \cap (G^0 \times G^0)$ then $u = s(u) = r(v) = v$
by property~(i) of Definition~\ref{df:gpoid}, and so $uv = us(u) = u$
by property~(ii). Therefore
$\tilde m\colon \mathcal{M}_{G^0}^{G^{(2)}} \to \mathcal{M}_{G^0}^G$
is smooth. It remains to prove that $\mathcal{M}_{G^0}^{G^{(2)}}$ is
diffeomorphic to $(\mathcal{M}_{G^0}^G)^{(2)}$. But this is clear on
examining the definitions: indeed,
\begin{eqnarray*}
\mathcal{M}_{G^0}^{G^{(2)}}
&=& G^{(2)} \times (0, \hbar_0] \uplus \mathcal{N}_{G^0}^{G^{(2)}}
\\
&=& \set{(g, h, \hbar): s(g) = r(h)}
    \uplus \set{(u,X_u,Y_u) : X_u, Y_u \in \mathcal{N}_{G^0}^G},
\end{eqnarray*}
whereas
\[
(\mathcal{M}_{G^0}^G)^{(2)}
 = \set{(g,\hbar_1; h,\hbar_2) : s(g) = r(h),\ \hbar_1 = \hbar_2}
   \uplus \set{(u,X_u; v,Y_v) : u = v}.
\]
For $G = G_M$, that boils down to
$TM \oplus TM \approx \mathcal{N}_M^{M\times M\times M}$. Note finally
the dimension count: the dimension of $(\mathcal{M}_{G^0}^G)^{(2)}$ is
$2(\dim G + 1) - \dim (G^0 + 1) = 2\dim G - \dim G^0 + 1$, which is
clearly the same as the dimension of $\mathcal{M}_{G^0}^{G^{(2)}}$.

This completes the proof of Proposition~\ref{pr:smooth-gpoid}.   

\section{Tangent groupoids and strict quantization}
\label{sc:three}

\subsection{The deformation conditions}

In our definition of quantization, we actually strengthen some of
Rieffel's requirements. Regard $T^*M$ as a Poisson manifold and
consider the classical $C^*$-algebra $A_0 := C_0(T^*M)$ of continuous
functions vanishing at infinity. We choose a dense subalgebra
$\mathcal{A}_0$ (there is considerable freedom in that, but, to fix
ideas, we think of the functions whose Fourier transform in the second
argument has compact support), and we search for a family of mappings
$Q_\hbar$ into noncommutative $C^*$-algebras $A_\hbar$ such that the
following relations hold for arbitrary functions in $\mathcal{A}_0$:
\begin{enumerate}
\item  the map $\hbar \mapsto \|Q_\hbar(f)\|$ is continuous on
       $[0,h_0)$ with $Q_0 = I$;
\label{it:one}
\item  $\lim_{\hbar\to 0} \|Q_\hbar(f_1) Q_\hbar(f_2) - Q_\hbar(f_2)
         Q_\hbar(f_1) - i\hbar Q_\hbar(\{f_1,f_2\}) \| = 0$;

\noindent
Those are Rieffel's \textit{strict quantization} conditions
(Question~23 of~\cite{RieffelProgram}), to which we add:

\item  the \textit{asymptotic morphism} condition $\lim_{\hbar\to 0}
        \|Q_\hbar(f_1)Q_\hbar(f_2) - Q_\hbar(f_1f_2)\| = 0$;
\item  the \textit{reality} condition $Q_\hbar(f^*) = Q_\hbar(f)^*$;
\label{it:four}
       and also
\item  the \textit{traciality} condition
       $\opname{Tr}[Q_\hbar(f_1) Q_\hbar(f_2)]
         = \int_{T^*M} f_1(q,p) f_2(q,p) \,\d\mu_\hbar(q,p)$;
\label{it:five}
\end{enumerate}
\noindent
where we use the same symbols $*$ (a bit overworked, admittedly) and
$\|\cdot\|$ for the adjoint and norm in every $C^*$-algebra.

Axioms (\ref{it:one}) to~(\ref{it:four}) are a slight variant of
Landsman's axioms. That the tangent groupoid construction provides an
answer to the twentieth query by Rieffel follows indeed from
Landsman's calculations in~\cite{Landsman}. There is no point in
repeating them here, and we limit ourselves to the necessary remarks
to fit them in the tangent groupoid framework. Axiom~(\ref{it:five})
is employed to further select a unique recipe.

In~\cite{RieffelProgram}, Rieffel introduces the important concept of
``flabbiness'': a deformation quantization is \textit{flabby} if it
contains the algebra of smooth functions of compact support on~$M$.
The constructions performed in this paper require in principle only
the existence of sprays, in order to use the tubular neighbourhood
theorem. However, in that case flabbiness is not guaranteed (see
below).

\subsection{The $C^*$-algebra of a groupoid}

The natural operation on functions of a groupoid is
\textit{convolution}:
\[
(a * b)(g) := \int_{\{hk=g\}} a(h)\, b(k)
 = \int_{\{h:r(h)=r(g)\}} a(h)\, b(h^{-1}g)
\]
but for this to make sense we need a measure to integrate with. We can
either define a family of measures on the fibres of the map~$r$,
$G^x := \set{g \in G : r(g) = x}$ for $x \in G^0$ (see the detailed
treatments given by Kastler~\cite{Kastler} and
Renault~\cite{Renault}) or we can finesse the issue by ensuring that
the integrand is always a $1$-density on each~$G^x$.

In the second approach, one uses half-densities, rather than
functions. We summarize it here, for completeness. Denote the typical
fibre of~$s$ by $G_y := \set{g \in G : s(g) = y}$ for $y \in G^0$.
Since $r$ and $s$ are submersions, the fibres $G^x$ and $G_y$ are
submanifolds of~$G$ of the same dimension, say~$k$. If $x = r(g)$ and
$y = s(g)$, then $\Lambda^k T_g G^x$ and $\Lambda^k T_g G_y$ are
lines. Let $\Omega_g^{1/2}$ be the set of maps
\[
\rho : \Lambda^k T_g G^x \otimes \Lambda^k T_g G_y  \to  \Cset
 \sepword{such that}
 \rho(t\alpha) = |t|^{1/2} \rho(\alpha) \quad\mbox{for $t \in \Rset$}.
\]
This is a (complex) line, and it forms the fibre at $g$ of a line
bundle $\Omega^{1/2} \to G$, called the ``half-density bundle''. Let
$C_c^\infty(G,\Omega^{1/2})$ be the space of smooth, compactly
supported sections of this bundle. For $a$, $b$ in this space, the
convolution formula makes sense and
$a * b \in C_c^\infty(G,\Omega^{1/2})$ also. The $C^*$-algebra of the
smooth groupoid $G \toto G^0$ is the algebra $C^*(G)$ obtained by
completing $C_c^\infty(G,\Omega^{1/2})$ in the norm
$\|a\| := \sup_{y\in G^0} \|\pi_y(a)\|$, where $\pi_y$ is the
representation of $C_c^\infty(G,\Omega^{1/2})$ on the Hilbert space
$L^2(G_y,\Omega^{1/2})$ of half-densities on the $s$-fibre $G_y$:
\[
\pi_y(a)\xi : g \mapsto \int_{G_y} a(h) \,\xi(h^{-1}g)
\]
where one notices that the integrand is a $1$-density on~$G_y$. If
$G = M \times M$, we get just the convolution of kernels:
\[
(a * b)(x,z) := \int_M a(x,y)\, b(y,z) \,\d y
\]
where $\d y$ denotes integration of a $1$-density on~$M$ parametrized
by~$y$.

This business of half-densities is very canonical and independent of
preassigned measures. However, in our case, if $M$ is an oriented
Riemannian manifold, we may use the volume form
$\,\d\nu(x) := \sqrt{g(x)}\,\d x$, where $g(x) := \det[g_{ij}(x)]$, to
replace half-densities by square-integrable functions on~$M \times M$
and on~$M$; thus $C^*(M \times M)$ is the completion of
$C_c^\infty(M \times M)$ acting as integral kernels on $L^2(M)$, so
that $C^*(M \times M) \simeq \mathcal{K}$, the $C^*$-algebra of
compact operators. If $G = TM$ is the tangent bundle, then $C^*(TM)$
is the completion of the convolution algebra
\[
(f_1 * f_2)(q, X)
 := \int_{T_qM} f_1(q, Y)\, f_2(q, X - Y) \,\sqrt{g(q)}\,\d Y
\]
where we may take $f_1(x,\cdot)$ and $f_2(x,\cdot)$ in
$C_c^\infty(T_xM)$. The \textit{Fourier transform}
\[
\mathcal{F} a(q,p) = \frac{1}{(2\pi)^n}
  \int_{T_qM} e^{-\ii pX} a(q,X) \,\sqrt{g(q)}\,\d X
\]
replaces convolution by the ordinary product on the total space $T^*M$
of the cotangent bundle. This extends to the isomorphism also called
$\mathcal{F} \colon C^*(TM) \to C_0(T^*M)$, with inverse:
\[
\mathcal{F}^{-1}b(q,X)
 = \int_{T_q^*M} e^{\ii pX} b(q,p) {\d^np \over \sqrt{g(q)}}.
\]

\subsection{The quantization and dequantization recipes}

Let $\gamma_{q,X}$ be the geodesic on~$M$ starting at~$q$ with
velocity~$X$, with an affine parameter~$s$, i.e.,
$\gamma_{q,tX}(s) \equiv \gamma_{q,X}(ts)$. Locally, we may write
\[
\begin{array}{rcl}
x & := & \gamma_{q,X}(s), \\
y & := & \gamma_{q,X}(-s),
\end{array}
 \sepword{with Jacobian matrix}  \pd{(x,y)}{(q,X)}(s).
\]

Then one has the change of variables formula:
\[
\int_{M\times M} F(x,y) \,\d\nu(x) \,\d\nu(y)
 = \int_M \int_{T_qM} F(\gamma_{q,X}(\half),\gamma_{q,X}(-\half))
     \,J(q,X;\half) \,\sqrt{g(q)} \,\d X \,\d\nu(q),
\]
where we introduce
\[
J(q,X;s) := s^{-n}
 \frac{\sqrt{\det g(\gamma_{q,X}(s))} \sqrt{\det g(\gamma_{q,X}(-s))}}
      {\det g(q)} \,\left| \pd{(x,y)}{(q,X)} \right|(s).
\]

This object can be computed from the equations of geodesic
deviation~\cite{Landsman}. The crucial estimate is
\begin{equation}
J(q,X,\half\hbar) = 1 + O(\hbar^2).
\label{eq:Jacobi-approx}
\end{equation}
In the computation of the Jacobian two types of Jacobi fields
intervene, let us call them $h$ and $\tilde h$, evaluated in different
points. The first is obtained through the variation in the coordinates
of the base point~$x$: in components,
\[
h_\nu^\mu(q,X,s) = \pd{x_1^\mu(q,X,s)}{q^\nu};
\]
the second one using the variation of the tangent vector,
\[
\tilde h_\nu^\mu(q,X,s) = \pd{x_1^\mu(q,X,s)}{X^\nu};
\]
and analogously for $x_2$ at~$-s$. To check the calculations, we may
work in normal coordinates, in which
\[
h_\nu^\mu = \delta_\nu^\mu(1 + O(\hbar^2)), \quad
\tilde h_\nu^\mu = \hbar\, \delta_\nu^\mu(1 + O(\hbar^2)),
\]
and moreover the metric has a simple expression.

Connes' quantization/dequantization maps are just the consequence of
applying the Gelfand--Na\u{\i}mark functor to the tangent groupoid
construction and are given simply by
\[
Q = \Phi_*^{-1} \circ J^{-1/2} \circ \mathcal{F}^{-1}
\]
and:
\[
Q^{-1} = \mathcal{F} \circ J^{1/2} \circ \Phi_* \,.
\]
Here $J^{-1/2}$, $J^{1/2}$ mean the corresponding multiplication
operators. It is possible to rewrite the formulae so $J$ only appears
once, but then we would lose property~(\ref{it:five}). So there is
actually a family of ``Moyal'' quantizations in the general case,
unless we demand traciality (as we do). This is a feature of the
nonflat case: note that $J = 1$ in the framework of our
Section~\ref{sc:one}. In the main development in~\cite{Landsman} $J$
actually only appears in the dequantization formula, but again
Landsman gives indication on how to modify the formulae to get
traciality. Those factors are analogous to the preexponential factors
that appear in the semiclassical expression for the path integral:
see, e.g., \cite[p.~95]{Schulman}. One can insert them in several ways
without altering the axioms, precisely on account
of~\eq{\ref{eq:Jacobi-approx}}.

A bit more explicitly, the previous formulae are
\begin{eqnarray}
k_a(x,y;\hbar)
&=& J^{-1/2}(q,X,\half\hbar) \,\mathcal{F}^{-1}a(q,X),
\nonumber \\
a(q,\xi)
&=& \mathcal{F}\bigl[ J^{1/2}(q,X,\half\hbar)\, k_a(x,y;\hbar) \bigr],
\label{eq:quantzn-maps}
\end{eqnarray}
where $(x,y)$ and $(q,X)$ are related by
$x = \gamma_{q,X}(\half\hbar)$, $y = \gamma_{q,X}(-\half\hbar)$.

The long but relatively straightforward verifications
in~\cite{Landsman} then show that we have defined an (obviously real)
preasymptotic morphism which moreover is a strict quantization from
$C_c^\infty(T^*M)$ to $\mathcal{K}(L^2(M))$. In addition, the tracial
property~(\ref{it:five}) is satisfied. In the Riemannian context we
have total tubular neighbourhoods, and it is clear that Fourier
transforms of smooth functions of compact support in $T^*M$ decay fast
enough that our formulae make sense for them; therefore the
quantization is \textit{flabby}.

We make a final comment on the uniqueness, or lack of it, of the
quantization considered here. Among our choices are those we make to
establish the isomorphism between $TM$ and $\mathcal{N}^\Delta$ on
which the whole construction hinges. Regarding only the differential
structure, those choices are parametrized by, say, the pair
$(\varphi_1 + \varphi_2, \varphi_1 - \varphi_2)
 \in \opname{End}(TM)\times GL(TM)$. However, the identification of
the normal bundle with the orthogonal bundle to the diagonal leads to
the equation $\varphi_1 + \varphi_2 = 0$, if we take the natural
metric on $M \times M$. To the same equation leads the reality
constraint~(\ref{it:four}), as a simple calculation
from~\eq{\ref{eq:quantzn-maps}} shows.

After so narrowing the freedom of parameter choice to $GL(TM)$,
semitraciality implies at least $\det(\varphi_1 - \varphi_2) = 1$.
Connes' condition~\eq{\ref{eq:bdry-tgt}} and our discussion in
Section~\ref{sc:one} strongly suggest to adopt the restriction
$\varphi_1 - \varphi_2 = 1$. Let us agree to call all the
quantizations (parametrized by $\opname{End}(TM)$) for which the last
equation holds quantizations \textit{of the Connes type}. Then we have
proved:

\begin{thm}
For any Riemannian manifold, the only\/ \textbf{real} quantization
rules of the Connes type are Moyal quantizations.   
\end{thm}

The classification of the various quantizations obtained by relaxing
the constraints imposed here deserves further investigation.

\section{Conclusions and Outlook}

The nicest aspect of Connes' construction is perhaps that the
continuity requirements in deformation quantization can be given a
canonical formulation in $C^*$-algebraic terms. In fact, we have a
continuous field of $C^*$-algebras in the sense of
Dixmier~\cite{Dixmier}, parametrized by $[0,\hbar_0]$. This
constitutes a \textit{strong} deformation of $C_0(T^*M)$
into~$\mathcal{K}$, i.e., the field of $C^*$-algebras is trivial
except at the origin.

There is a short exact sequence of $C^*$-algebras
\[
0 \longto C_0(0,\hbar_0] \otimes \mathcal{K} \longto C^*(G_M)
  \longto^\sigma C_0(T^*M) \longto 0
\]
that yields isomorphisms in $K$-theory:
\[
K_j(C^*(G_M)) \longto^{\sigma_*} K_j(C_0(T^*M)) = K^j(T^*M)
 \qquad  (j = 0,1).
\]
This corresponds to Proposition~II.5.5 of Connes~\cite{Book} and the
claims immediately thereafter. It can be seen as follows: given a
smooth groupoid $G = G_1 \uplus G_2$ which is a disjoint union of two
smooth groupoids with $G_1$ open and $G_2$ closed in~$G$, there is a
short exact sequence of $C^*$-algebras
\[
0 \longto C^*(G_1) \longto C^*(G) \longto^\sigma C^*(G_2) \longto 0
\]
where $\sigma$ is the homomorphism of restriction from $C_c^\infty(G)$
to $C_c^\infty(G_2)$: it is enough to notice that $\sigma$ is
continuous for the $C^*$-norms because one takes the supremum of
$\|\pi_u(a)\|$ over the closed subset $u \in G^{(2)}_2$, and it is
clear that $\ker\sigma \simeq C^*(G_1)$. Since
$C^*(M \times M) = \mathcal{K}$, the $C^*$-algebra $C^*(G_1)$,
obtained by completing the algebraic tensor product
$C_c^\infty(0,\hbar_0] \odot C_c^\infty(M \times M)$, is
$C_0(0,\hbar_0] \otimes \mathcal{K}$, which is contractible via the
homotopy $\alpha_t(f \otimes A) := f(t\,\cdot) \otimes A$ for
$f \in C_0(0,\hbar_0]$, $0 \leq t \leq 1$; in particular,
$K_j(C_0(0,\hbar_0] \otimes \mathcal{K}) = 0$. At this point, we
appeal to the six-term cyclic exact sequence in $K$-theory of
$C^*$-algebras~\cite{Wegge}; the two trivial groups break the circuit
and leave the two advertised isomorphisms.

The restriction of elements of $C^*(G)$ to the outer boundary
$M \times M \times \{\hbar_0\}$ gives a homomorphism
\[
\rho :
 C^*(G_M) \to C^*(M \times M \times \{\hbar_0\}) \simeq \mathcal{K},
\]
and in $K$-theory this yields a homomorphism
$\rho_*\colon K_0(C^*(G_M)) \to K_0(\mathcal{K}) = \Zset$. Finally, we
have the composition
$\rho_*(\sigma_*)^{-1}\colon K^0(T^*M) \to \Zset$, which is just the
Atiyah--Singer analytical index map. (To see that, one must identify
suitable pseudodifferential operators with the kernels in $C^*(G_1)$:
see~\cite{MonthubertAmiPierrot,MonthubertII}; of course, general
results in the theory of pseudodifferential operators~\cite{Petersen}
ensure that the index does not depend on the particular quantization
recipe chosen.)

To summarize, we found Connes' tangent groupoid construction to
address Rieffel's questions concerning the r\^ole of a Riemannian
metric in a strict (deformation) quantization. The presence of a
Riemannian metric on a manifold~$M$ allows one to identify a geodesic
spray to construct a total tubular neighbourhood; to choose a standard
representative for the normal bundle, to have a natural isomorphism
between the normal and the tangent bundle; to use a canonical measure
on~$M$; and to define naturally the fibrewise Fourier transform.

Work on extensions of some of the methods in this paper to more
general contexts is in progress.

\begin{ack}
Soon after a previous version of this manuscript was issued, we
received a comment from Marc Rieffel clarifying the difference between
strict quantization and strict deformation quantization. This helped
to make some of the statements in this paper more precise; we are very
grateful to him for that.

JMGB thanks the Departamento de F\'{\i}sica Te\'orica of the
Universidad de Zaragoza for a customarily warm hospitality; he and JCV
acknowledge support from the Vicerrector\'{\i}a de Investigaci\'on of
the University of Costa Rica. Support of Spanish DGICYT (PB93--0582)
and DGES (PB96--0717) is also acknowledged.
\end{ack}

\end{document}